\documentclass{mypnastwo} 
\usepackage[T1]{fontenc}
\usepackage[utf8]{inputenc}
\usepackage{amssymb,amsmath,amsfonts}
\usepackage[noadjust]{cite}
\usepackage{array}
\usepackage{booktabs}
\usepackage{mdwtab}
\usepackage{mathtools}
\usepackage{mathptmx}
\usepackage{hyphenat}
\usepackage{enumitem}
\usepackage{ifpdf}
\usepackage{amsthm}
\usepackage{sansmath}
\usepackage{pnastwof}
\ifpdf
  \usepackage[pdftex]{graphicx}
  \usepackage[bookmarks=true, bookmarksopen=true,%
    bookmarksdepth=3,bookmarksopenlevel=2,%
    colorlinks=true,%
    linkcolor=blue,%
    citecolor=blue,%
    filecolor=blue,%
    menucolor=blue,%
    urlcolor=blue]{hyperref}
  \hypersetup{pdftitle={A positive basis for surface skein algebras}}
  \hypersetup{pdfauthor={Dylan P. Thurston}}
\else
  \usepackage[dvips]{graphicx}
  \DeclareGraphicsExtensions{.eps,.ps,.mps}
  \usepackage[draft]{hyperref}
\fi
\usepackage{url}

\makeatletter
\newcommand\mi@kern[1]{%
  \settowidth\@tempdima{$\mi@obj^{#1}$}
  \kern-\@tempdima
  #1
  \settowidth\@tempdima{$\mi@obj$}
  \kern\@tempdima
}

\newtoks\mi@toksp
\newtoks\mi@toksb
\DeclareRobustCommand{\manyindices}[5]{
  \def\mi@obj{#5}
  \mi@toksp\expandafter{\mi@kern{#2}}
  \mi@toksb\expandafter{\mi@kern{#1}}
  \@mathmeasure4\textstyle{#5_{#1}^{#2}}
  \@mathmeasure6\textstyle{#5_{#3}^{#4}}
  \dimen0-\wd6 \advance\dimen0\wd4
  \@mathmeasure8\textstyle{\hphantom{{}_{#1}^{#2}}#5^{\the\mi@toksp#4}_{\the\mi@toksb#3}}
  \hbox to \dimen0{}{\kern-\dimen0\box8}
}
\makeatother 

\newcommand{\margin}[1]{%
  \marginpar[\raggedleft\footnotesize{#1}]{\raggedright\footnotesize{#1}}}

\newcommand{\move}[2][\longrightarrow]{\stackrel{\text{#2}}{#1}}

\newread\testin

\def\mathcenter#1{%
  \vcenter{\hbox{$#1$}}%
}

\def\mfig#1{
        \mathcenter{\includegraphics{#1}}
}

\def\mfigb#1{
        \mathcenter{\includegraphics[trim=-1 -1 -1 -1]{#1}}
}

\newcommand{\RR}{\mathbb R}

\newcommand{\ZZ}{\mathbb Z}
\newcommand{\QQ}{\mathbb Q}

\renewcommand{\sl}{\mathfrak{sl}}

\newcommand{\co}{\colon}

\newcommand{\abs}[1]{{\lvert #1 \rvert}}

\newcommand{\bdy}{\partial}

\DeclareMathOperator{\tr}{tr}

\DeclareMathOperator{\lcm}{lcm}

\DeclareMathOperator{\SL}{\mathrm{SL}}

\newtheoremstyle{plain-pnas}%
  {\abovedisplayskip}
  {\belowdisplayskip}
  {\itshape}
  {}
  {\theoremfont}
  {.}
  { }
  {#1 #2\thmnote{ (#3)}}
\newtheoremstyle{def-pnas}%
  {\abovedisplayskip}
  {\belowdisplayskip}
  {}
  {}
  {\theoremfont}
  {.}
  { }
  {#1 #2\thmnote{ (#3)}}
\theoremstyle{plain-pnas}
\newtheorem{theorem}{Theorem}
\newtheorem{proposition}{Proposition}[section]
\newtheorem{citethm}[proposition]{Theorem}
\newtheorem{lemma}[proposition]{Lemma}
\newtheorem{corollary}[proposition]{Corollary}

\newtheorem{conjecture}[proposition]{Conjecture}

\newtheorem{definition}[proposition]{Definition}

\newtheorem{question}[proposition]{Question}
\newtheorem{problem}[proposition]{Problem}

\newtheorem{example}[proposition]{Example}
\newtheorem{remark}[proposition]{Remark}

\newcommand{\Surf}{\mathbf{S}}
\newcommand{\Mark}{\mathbf{M}}
\newcommand{\SM}{(\Surf,\Mark)}

\newcommand{\Diag}{\mathbf{D}}
\newcommand{\Curve}{\mathbf{C}}

\newcommand{\Skein}{\mathrm{Sk}}
\newcommand{\Sk}{\Skein}
\newcommand{\Basis}{\mathbf{B}}
\DeclareMathOperator{\Bangle}{Bang}
\newcommand{\Bang}{\Bangle}
\DeclareMathOperator{\Band}{Band}
\DeclareMathOperator{\Bracelet}{Brac}
\newcommand{\Brac}{\Bracelet}
\DeclareMathOperator{\Cross}{Cross}

\renewcommand{\t}[1]{\widetilde{#1}}

\providecommand{\noopsort}[1]{}

\renewcommand{\margin}[1]{}

\graphicspath{{draws/}}

\def\theurl{}
\copyrightyear{2013}
\issuedate{Issue Date}
\volume{Volume}
\issuenumber{Issue Number}
\footlineauthor{Thurston}

\begin{document}
\title{A positive basis for surface skein algebras}


\author{Dylan P. Thurston\affil{1}{Department of Mathematics, Indiana
    University, Bloomington, IN 47405}}
\contributor{Submitted to Proceedings of the National Academy of Sciences\\
Dedicated to the memory of William P. Thurston, 1946--2012}
\date{\today}

\maketitle

\begin{article}
\begin{abstract}
  We show that the twisted $\SL_2$ skein algebra of a surface
  has a natural basis (the
  bracelets basis) that is positive, in the sense that the structure
  constants for multiplication are positive integers.
\end{abstract}

\keywords{skein algebra | character variety | canonical basis | positive basis}

\setcounter{tocdepth}{1}

The Jones polynomial of knots is one of the simplest and most powerful
knot invariants, at the center of many recent advances in topology; it
is a polynomial in a parameter $q$. The skein algebra of a surface is
a natural generalization of the Jones polynomial to knots that live in
a thickened surface. In this paper, we propose a new basis for the
skein algebra. This basis has positivity properties when $q$ is set to
1, and conjecturally for general values of $q$ as well. This is part
of a more general conjecture for cluster algebras, and suggests the
existence of well-behaved higher-dimensional structures.

\section{Introduction}
\label{sec:intro}

\dropcap{F}or a compact oriented surface $\Sigma$ (possibly with boundary), the
\emph{Kauffman bracket skein algebra}, denoted $\Sk_q(\Sigma)$, is the
$\ZZ[q^{\pm 1}]$-module
spanned by framed links in $\Sigma \times [0,1]$ modulo the local
relations
\begin{align}
  \label{eq:kauffman}
  \left\langle \mfig{curves-3} \right\rangle
   &= 
  q \left\langle \mfig{curves-2} \right\rangle +
  q^{-1} \left\langle \mfig{curves-1} \right\rangle\\
  \label{eq:unknot}
  \left\langle \mfig{curves-10} \right\rangle &= -q^2 -q^{-2}.
\end{align}
Vertical stacking of links makes $\Sk_q(\Sigma)$ into an algebra: to
form $\langle D_1\rangle \cdot \langle D_2\rangle$, superimpose $D_1$
onto $D_2$, making $D_1$ cross over $D_2$.

This skein algebra was first defined by
Przytycki~\cite{Przytycki91:SkeinMod} and
Turaev~\cite{Turaev89:AlgebraLoops} as an extension of the Jones
polynomial of knots in $S^3$ to knots in a surface cross an interval.
When specialized to $q=\pm1$, 
we no longer need to record crossing information.  For $q=-1$, we
essentially get the
algebra of functions on the $\SL_2(\RR)$ character variety of~$\Sigma$
\cite{Bullock97:RingsSL2,PS00:SkeinSL2,MS13:CharSL2}.
A choice of a
spin structure gives an isomorphism between $\Sk_q(\Sigma)$ and
$\Sk_{-q}(\Sigma)$~\cite{Barrett99:SkeinSpin}.  More naturally, $\Sk_1(\Sigma)$
can be thought of as the algebra of functions on the twisted
$\SL_2(\RR)$ character variety.

\begin{definition}
  A \emph{twisted $\SL_2(\RR)$ representation} of a surface~$\Sigma$
  is a representation of $\pi_1(UT\Sigma)$, the fundamental group of
  the unit tangent bundle of~$\Sigma$, into $\SL_2(\RR)$, with the
  property that rotation by $2\pi$ acts by $-1 \in
  \SL_2(\RR)$.  The \emph{twisted $\SL_2(\RR)$ character variety} is
  the algebro-geometric quotient of twisted $\SL_2(\RR)$
  representations by conjugation.
\end{definition}

A hyperbolic structure on $\Sigma$ gives a canonical
twisted $\SL_2(\RR)$ representation.  See, e.g.,
\cite[Prop.~10]{BW11:QuantTraces}.

In this paper, we are mainly interested in $\Sk_1(\Sigma)$, henceforth
denoted $\Sk(\Sigma)$.
Our main result is that it has a \emph{positive
basis}.
\begin{definition}\label{def:positive}
  For an algebra $A$ over $\ZZ$ (free as a $\ZZ$-module), a basis
  $\{x_i\}$ is \emph{positive} if
  \[
  x_i \cdot x_j = \sum_k m_{ij}^k x_k
  \]
  where $m_{ij}^k \ge 0$.
\end{definition}

\begin{figure}
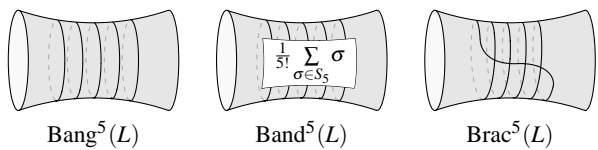

  \vskip-\abovedisplayskip
  \begin{align*}
    &\underset{\displaystyle\Bang^5(L)}{\mfigb{annulus-0}}
    &&\underset{\displaystyle\Band^5(L)}{\mfigb{annulus-1}}
    &&\underset{\displaystyle\Brac^5(L)}{\mfigb{annulus-2}}
  \end{align*}
  \caption{Examples of bangle, band, and bracelet operations applied
    to the core loop of an annulus.  The bangle has parallel copies,
    the band averages over all ways of joining, and the bracelet
    wraps multiple times.  \label{fig:three-bases}}
\end{figure}

We will show that the \emph{bracelets basis} (Def.~\ref{def:bases}) of
the skein algebra is
positive. This basis
is \emph{not} made of crossingless diagrams. In
Fig.~\ref{fig:three-bases}, instead of \emph{bangles}
we use \emph{bracelets}.

\begin{theorem}\label{thm:positive-1}
  The bracelets basis is a natural
  positive basis for $\Sk(\Sigma)$.
\end{theorem}

The basis is \emph{natural} in the sense that it is invariant under the
mapping class group (automorphisms of the surface).  Although a spin
structure gives an isomorphism between
$\Sk_1(\Sigma)$ and $\Sk_{-1}(\Sigma)$ as algebras, it is unlikely
that $\Sk_{-1}(\Sigma)$ has a natural positive
basis, as $\Sigma$ generally does not have a canonical spin
structure.

We work with a mild extension of the skein algebra, to
include marked points and arcs with endpoints on the marked points.
In a sequel to this paper, we will extend to
connect this paper to cluster algebras, the
original motivation for this work.  Specifically, we will construct a
positive, natural basis for an algebra between the cluster
algebra and the upper cluster algebra of a marked surface, building on
the connections between surfaces and cluster algebras \cite{GSV05:ClusterAlgWPForms,FG09:ClusterEnsemb,FST08:ClusterI}.

Positivity was first
conjectured by Fock and Goncharov in their ground-breaking paper \cite[Section
12]{FG06:ModuliLocalSystems}.
The bracelets basis was considered 
by Musiker,
Schiffler, and Williams~\cite{MSW13:BasesSurfaces}, who
proved a weaker form of positivity.
This weaker positivity and explicit combinatorial formulas have been
well-studied
\cite{MS10:ClusterExpansion,MSW11:Positivity,MW13:Matrix,ST09:CustSurfI}.



In fact, we will prove a stronger theorem.

\begin{theorem}\label{thm:positive-expansion}
  For any diagram~$D$ on~$\Sigma$, the expansion
  of~$\langle D \rangle$ in the
  bracelets basis is sign-coherent.  If $D$ has no null-homotopic
  components or nugatory crossings, then the expansion is positive.
\end{theorem}

Here, \emph{sign-coherent} means that either all terms are positive or
all are negative.  A \emph{nugatory crossing} is a crossing that cuts
off a null-homotopic loop.

The proof proceeds by carefully picking a
crossing to resolve by Eq.~\ref{eq:kauffman}, being careful
to avoid ever introducing a negative sign.

\subsection{Extensions and future work}
\label{sec:extensions}

Theorem~\ref{thm:positive-1} suggests many possible extensions.
First, can anything be said when $q \ne 1$?

\begin{conjecture}\label{conj:quantum-skein}
  The bracelets basis can be lifted to a positive basis for the
  quantum skein algebra $\Sk_q(\Sigma)$.
\end{conjecture}
See Conj.~\ref{conj:bracelets-q-positive} for a more precise version.
This conjecture was essentially made by Fock and Goncharov
\cite[Conj.~12.4]{FG06:ModuliLocalSystems}.
The techniques in this paper will not work for the $q$-deformation, as
there is no obvious
analogue of Theorem~\ref{thm:positive-expansion} for the quantum skein
algebra.

Fig.~\ref{fig:three-bases} suggests another natural basis for $\Sk_q(\Sigma)$,
the \emph{bands basis}.

\begin{question}
  When is the bands basis positive?
\end{question}

See Conjecture~\ref{conj:bands-positive} for a possible answer.

Finally, the existence of a positive basis suggests the presence of a
``nice'' categorification, where product becomes a monoidal tensor
product and sum becomes direct sum, or possibly a composition series.
We leave the precise formulation vague.

\begin{question}
  Is there a monoidal categorification of $\Sk_1(\Sigma)$ or
  $\Sk_q(\Sigma)$ that makes positivity of the bracelets basis or
  bands basis manifest?
\end{question}



\section{Preliminaries on surfaces and curves}
\label{sec:curves}

We first collect some basic facts on the topology of curves in surfaces.

\begin{definition}
  A \emph{marked surface} $\Sigma = \SM$ is a pair of a surface $\Surf$,
  possibly with boundary, and a finite set $\Mark \subset \Surf$ of
  marked points.
  Marked points in the interior of $\Surf$ are called
  \emph{punctures}.
  When we think of $\Sigma$ itself as a topological space, we mean
  $\Surf \setminus \Mark$.  For convenience, we will assume that
  $\Surf$ is connected.
\end{definition}

Simple surfaces with few or no
marked points are not excluded.

\begin{definition}
  A \emph{(curve) diagram}~$D = (X,\phi)$ on $\Sigma$ is 1-manifold
  with boundary, $X$, and
  an immersion $\phi\co(X, \bdy X) \to (\Surf, \Mark)$, by which we mean an
  immersion of
  $X$ in $\Surf$ so that each boundary point of~$X$ maps to a marked
  point and no point in the interior of $X$ maps to a marked point.
  We also require that $D$ has only simple transverse crossings, and
  that no point in the interior of~$X$ maps to $\bdy\Surf$.
  $D$ is \emph{connected} if $X$ is connected, $D$ is an \emph{arc}
  if $X$ is an interval, and $D$ is a \emph{loop} if
  $X$ is a circle.
\end{definition}

There are several equivalence relations on diagrams.

\begin{figure}
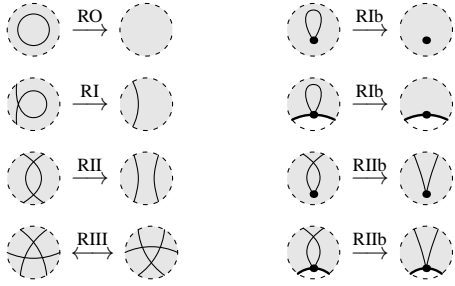

  \vskip-\abovedisplayskip
  \begin{align*}
    \mfigb{curves-10} &\move{RO} \mfigb{curves-12}&
      \mfigb{curves-44} &\move{RIb} \mfigb{curves-45}\\
    \mfigb{curves-60} &\move{RI} \mfigb{curves-63}&
      \mfigb{curves-42} &\move{RIb} \mfigb{curves-43}\\
    \mfigb{curves-70} &\move{RII} \mfigb{curves-72}&
      \mfigb{curves-80} &\move{RIIb} \mfigb{curves-82}\\
    \mfigb{curves-90} &\move[\longleftrightarrow]{RIII} \mfigb{curves-91}&
      \mfigb{curves-83} &\move{RIIb} \mfigb{curves-85}
  \end{align*}
  \caption{The Reidemeister moves. The diagrams show a local portion
    of the surface. The surface is shaded gray; white regions are not
    in the surface. In moves RIb and RIIb, there may be more arcs
    ending at the marked point, not intersecting the displayed
    arcs.\label{fig:reidemeister}}
\end{figure}

\begin{definition}
  A \emph{Reidemeister move} is one of the moves in
  Fig.~\ref{fig:reidemeister} (in either direction).  Note that this
  includes some moves that change the number of components.
  A Reidemeister \emph{reduction} is a Reidemeister move that reduces
  or keeps constant the number of intersections and components, i.e.,
  all moves from left to right, and RIII in either direction.
  A \emph{strict reduction} is any reduction other
  than RIII\@.
\end{definition}

\begin{definition}
  Two diagrams $D_1, D_2$ are \emph{ambient isotopic} if they can be
  related by an isotopy of~$\Sigma$, or equivalently if there is an
  isotopy of $D_1$ to $D_2$
  that does not change any of the crossings.  We always consider
  diagrams up to ambient isotopy; the set of diagrams up
  to ambient isotopy is denoted $\Diag(\Sigma)$.

  $D_1$ and $D_2$ are \emph{regular isotopic}
  if they can be connected by a path within the space of
  immersions, dropping the condition on crossings but keeping
  the interior of the diagram disjoint from $\bdy \Surf \cup \Mark$.
  Equivalently, $D_1$ can be connected to $D_2$ using moves RII\@, RIIb,
  and RIII\@.
  The set $\Diag(\Sigma)$ modulo regular isotopy is
  denoted $\Curve(\Sigma)$.

  $D_1$ and $D_2$ are \emph{homotopic}, written $D_1 \sim D_2$, if they can be connected within
  the space
  of all continuous maps.  
  Equivalently, $D_1$ can be connected to $D_2$ using moves RI\@, RII\@, RIIb,
  and RIII\@.

  $D_1$ \emph{reduces} to $D_2$ if $D_1$ can be turned in to $D_2$ by
  a chain of zero or more Reidemeister reductions, of any of the types.
\end{definition}

There is a product on diagrams modulo regular isotopy.

\begin{definition}
  For $C_1, C_2 \in \Curve(\Sigma)$, their \emph{product} $C_1 C_2$ is
  obtained by taking diagrams $D_1, D_2$ so that $C_i =
  [D_i]$ and $D_1$ intersects $D_2$ transversally, and defining
  \[
  C_1 C_2 = [D_1 \cup D_2].
  \]
\end{definition}

\begin{lemma}\label{lem:curves-prod}
  The product above is well-defined, and makes $\Curve(\Sigma)$ into a
  commutative, associative monoid with unit the empty diagram.
\end{lemma}

\begin{proof}
  Standard.
\end{proof}

We next give conditions for a diagram to have
minimal self-intersections.

\begin{definition}\label{def:chains-gons}
  For a diagram~$D = (X, \phi)$, a \emph{segment} of~$D$ is an
  oriented subinterval of~$X$ whose endpoints are either endpoints of $X$ or
  preimages of crossings.  If two segments $S_1$ and $S_2$ meet at a
  marked point~$p$, then a \emph{turn} from $S_1$ to~$S_2$ at~$p$ is a
  homotopy class of arcs from $\phi(S_1)$ to $\phi(S_2)$ inside $N \setminus p$,
  where $N$ is a small neighborhood of~$p$.  (If $p$ is not a
  puncture, there is only one such turn.)

  A \emph{$k$-chain} of~$D$ is a sequence of
  $k$~segments $(S_i)_{i=1}^k$ so that the end of~$S_i$
  has the same image as
  the start of~$S_{i+1}$, with $i$ interpreted modulo~$k$,
  and so that the~$S_i$ are disjoint subsets of~$X$.  Furthermore
  choose a turn from $S_i$ to $S_{i+1}$ whenever they meet at a marked
  point.
  As a special case, a \emph{$0$-chain} of~$D$ is a loop component
  of~$D$.

  For any
  $k$-chain~$H$, there is an associated loop $H^\circ$, obtained by
  smoothing out the corners at the endpoints of the segments. At
  marked points, follow the chosen turn without
  nugatory crossings:
  \begin{equation}
    \label{eq:turn}
    \mfigb{curves-120} \mapsto \mfigb{curves-121}
  \end{equation}
  An \emph{embedded $k$-gon} is a $k$-chain~$H$ so that $H^\circ$ bounds
  a disk.  A \emph{singular $k$-gon} or just \emph{$k$-gon} is a
  $k$-chain~$H$ so that $H^\circ$ is null-homotopic.  Finally, a
  \emph{weak segment} is an immersed subinterval~$S$ of~$X$ whose
  endpoints are crossings of~$D$ or endpoints of~$X$.  (That is, $S$
  may wrap more than once around a loop component of~$D$.)  A \emph{weak
  $k$-chain} is like a $k$-chain, but using weak segments and dropping
  the requirement that the $S_i$ be disjoint, and a \emph{weak
    $k$-gon} is a weak $k$-chain~$H$ so that $H^\circ$ is
  null-homotopic.  (For $H$ a weak $k$-chain, $H^\circ$ is defined
  analogously to the above definition, possibly running multiple times
  over some portions of~$D$.) $0$-gons,
  $1$-gons, and $2$-gons are also called disks, monogons, and
  bigons. See Fig.~\ref{fig:bigons}.

  \begin{figure}
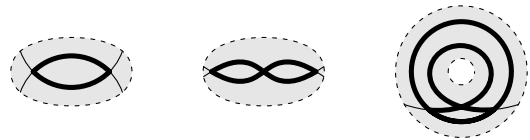

    \vskip-\abovedisplayskip
    \begin{align*}
      &\mfig{curves-110}
      &&\mfig{curves-111}
      &&\mfig{curves-112}
    \end{align*}
    \caption{Embedded, singular, and weak bigons. Two portions of the
      weak bigon run over the same part of the loop.\label{fig:bigons}}
  \end{figure}

  To relate chains to the fundamental group, note that a
  representative of
  $\alpha \in \pi_1(\Sigma, x)$ can be viewed as a $1$-chain, so
  $\alpha^\circ$ is a loop on $\Sigma$.  Conversely, the
  \emph{holonomy} of a loop~$L$ on~$\Sigma$ is the corresponding
  element of $\pi_1(\Sigma, x)$, where we connect~$L$ to the
  basepoint~$x$ by a specified path.
  (The term ``holonomy'' comes from thinking about the canonical
  $\pi_1(\Sigma)$-bundle over~$\Sigma$.)
\end{definition}

For example, the left hand side of
\begin{itemize}
\item RO has an embedded disk,
\item RI and RIb have an embedded monogon,
\item RII and RIIb have an embedded bigon, and
\item RIII has an embedded triangle.
\end{itemize}

\begin{definition}
  In a diagram~$D$,
  a \emph{trivial component} is a
  null-homotopic component, necessarily either a loop or an arc from a
  marked point to itself.  $D$ is \emph{simple} if it has no
  crossings and no trivial components.
  $D$ is \emph{taut} if it has a
  minimum number of self-intersections in its homotopy class and has
  no trivial components.
\end{definition}

\begin{lemma}\label{lem:taut-union}
  If $D = D_1 \cup D_2$ is a taut diagram, then $D_1$ and $D_2$ are
  also taut diagrams.
\end{lemma}

\begin{proof}
  This follows, for instance, from the fact that a diagram is taut iff
  it is length-minimizing with respect to some metric \cite{Shepard91:TopologyShortest,NC01:ShortestGeodesics}.
\end{proof}

Diagrams can be
monotonically simplified:

\begin{citethm}[Hass-Scott~\cite{HS94:ShorteningCurves}]
  \label{thm:curve-shorten}
  An arbitrary diagram can be turned into a taut diagram by a
  sequence of Reidemeister reductions.
\end{citethm}

(Hass and Scott proved a version of this theorem without arcs or marked
points, but the techniques extend immediately.)

\begin{definition}
  For $D$ a diagram, the set of crossings of~$D$ is denoted $\Cross(D)$.
  For $x\in\Cross(D)$, a \emph{resolution}
  of~$x$ is one of the two diagrams obtained by replacing a
  neighborhood of~$x$ by a local picture without crossings.
  A \emph{connected} resolution is one that does not increase the
  number of components of~$D$.
\end{definition}

\begin{lemma}\label{lem:cross-correspond}
  Let $S$ be a set of Reidemeister reductions containing RII\@.
  If $D_1$ reduces to $D_2$ by a set of reductions from~$S$, then there is an
  injection~$m: \Cross(D_2) \to \Cross(D_1)$ so that for each~$x \in
  \Cross(D_2)$, the
  resolutions of
  $m(x)$ are related by moves in~$S$ to resolutions of~$x$.
\end{lemma}

\begin{proof}
  It suffices to consider the case when $D_1$ and $D_2$ are related by
  a single Reidemeister reduction.  When $x \in \Cross(D_2)$ is outside of the
  reducing region, $m(x)$ is the same crossing as an
  element of $\Cross(D_1)$.  The only case where $x$ can be inside the
  reducing region is for RIII\@, where $m$ is the correspondence
  \[
  \mfigb{curves-95} \longleftrightarrow \mfigb{curves-94}.
  \]
  In this RIII case, the resolution at~$x$ and the resolution at $m(x)$ are
  related by two RII moves.  Otherwise, the resolutions are related by
  the same Reidemeister move that relates $D_1$ and $D_2$.
\end{proof}

\begin{lemma}\label{lem:no-gons}
  A taut diagram has no singular disks or monogons, and the only
  singular bigons are bigons between
  isotopic arcs.
\end{lemma}

\begin{proof}
  Let $D$ be a taut diagram.  It follows by
  definition that $D$ has no singular disk.  If $D$ has a
  singular monogon, then the connected resolution of the crossing of
  the monogon yields a
  diagram~$D'$ homotopic to~$D$ but with fewer crossings.  Similarly, if
  $D$ has a singular bigon involving at least one crossing, resolving the 
  crossing(s) yields a simpler diagram~$D' \sim D$.
\end{proof}

There is also a partial converse to the above lemma.

\begin{citethm}
\label{thm:gon-existence}
  If a diagram~$D$ is not taut, then it has a singular
  disk, singular monogon, or weak bigon.
\end{citethm}

\begin{proof}
  This was first proved by Hass and Scott in 1985
  \cite[Theorem~3.5]{HS85:IntersectionsCurves}.  We give a short
  proof using their later curve-shortening
  techniques~\cite{HS94:ShorteningCurves}, as we will
  use the technique for Lemma~\ref{lem:strong-weak-positive}
  below.  By Theorem~\ref{thm:curve-shorten}, there
  is a sequence of Reidemeister reductions
  \[
  D = D_0 \to D_1 \to \dots \to D_n
  \]
  where $D_n$ is taut.  Since $D_0$ is not taut, at least one of these
  reductions, say $D_k \to D_{k+1}$, is strictly reducing.  If $D_k
  \to D_{k+1}$ is:
  \begin{description}
  \item[RO:] $D_k$ has an embedded disk
  \item[RI:] $D_k$ has an embedded monogon
  \item[RII or RIIb:] $D_k$ has an embedded bigon.
  \end{description}
  In each case, we follow the disk, monogon, or bigon backwards
  from~$D_k$ to $D_0$, using the lemma below.
\end{proof}

\begin{lemma}\label{lem:gons-exist}
  If a diagram $D_1$ reduces to $D_2$ and $D_2$ has a singular disk,
  singular monogon, or weak bigon, then $D_1$ has one as well.
\end{lemma}

\begin{proof}
  We use the map $m\co \Cross(D_2) \to \Cross(D_1)$ from
  Lemma \ref{lem:cross-correspond}.
  \begin{itemize}
  \item If $D_2$ has a singular disk, the homotopy type of this
    component is the same in~$D_1$.
  \item If $D_2$ has a singular monogon with corner $x$, then $m(x)$
    is the corner of a singular monogon for~$D_1$.
  \item If $D_2$ has a weak bigon with corners at $x$ and $y$, then
    $D_1$ has a weak bigon with corners at $m(x)$ and $m(y)$.\qedhere
  \end{itemize}
\end{proof}

\begin{remark}\label{rem:weak-bigons}
  Theorem~\ref{thm:gon-existence} is false if we replace ``weak bigon''
  with ``singular bigon'' in the statement.  Let $L$ be a simple loop
  and $\Brac^3(L)$ the loop that wraps $3$ times around~$L$. Then
  there is a non-taut embedding
  of $L \cup \Brac^3(L)$ with a weak bigon but no singular bigon
  \cite{HS85:IntersectionsCurves}. Fig.~\ref{fig:bigons} on the right shows an
  example of problematic local behavior just after a singular bigon
  has become a weak bigon.
  On the other hand, a taut diagram may have a weak bigon; for instance,
  $\Brac^2(L)$ always has a weak bigon.  So
  Theorem~\ref{thm:gon-existence} does not give necessary and
  sufficient conditions for a diagram to be taut.
\end{remark}

Despite Remark~\ref{rem:weak-bigons}, there is a converse for
connected diagrams.

\begin{citethm}[{\cite[Lemma 4.1 and Theorem
    4.2]{HS85:IntersectionsCurves}}]\label{thm:gon-exist-2}
  If a connected diagram is not taut, then it has a singular
  disk, monogon, or bigon.
\end{citethm}

\begin{remark}
  It is possible to give a curve-shortening proof of
  Theorem~\ref{thm:gon-exist-2} along the lines of the proof of
  Theorem~\ref{thm:gon-existence} above.  The bad behavior in
  Fig.~\ref{fig:bigons} cannot occur when the diagram is
  connected.
\end{remark}

\emph{Bracelets}, loops that wrap
multiple times, play a key role in our
construction.
\begin{proposition}\label{prop:powers-intersect}
  If $A$ and $B$ are two taut diagrams with holonomy $\alpha$ and
  $\beta$,
  $B$ is simple, and $\alpha = \beta^k$,
  then $A$ has
  $k-1$ self-intersections.  We can label these
  self-intersections $x_i$, $1 \le i < k$, so the holonomy of the
  two components of the disconnected resolution of~$A$ at $x_i$ are
  conjugate to
  $\beta^i$ and $\beta^{k-i}$.
\end{proposition}

\begin{proof}
  See \cite[Lemma~1.9]{HS85:IntersectionsCurves} and
  \cite[Theorem~2.1, Case~3]{HS94:ShorteningCurves}.
\end{proof}

Theorem~\ref{thm:positive-expansion} says that taut diagrams have
positive expansions in
the bracelets basis.  There is a more general class of diagrams with
positive expansions.

\begin{definition}
  A diagram is \emph{weakly positive} if it is regular isotopic to a taut
  diagram.
  A diagram is \emph{strongly positive} if it has no singular
  monogons or singular disks.
\end{definition}

\begin{lemma}\label{lem:strong-weak-positive}
  A strongly positive diagram is weakly positive.
\end{lemma}

\begin{proof}
  Let $D$ be strongly positive.  By Theorem~\ref{thm:curve-shorten},
  $D$ can be reduced to a taut
  diagram~$D'$ using reductions.  By Lemma~\ref{lem:gons-exist}, if at
  any step along the way we use a RO move, then the singular disk can
  be followed back to give a singular disk for~$D$; but $D$ is
  strongly positive, so this cannot happen.  Similarly for singular monogons
  and RI or RIb moves.  Thus we only use moves RII\@, RIIb, and RIII
  in the reduction from $D$ to~$D'$, as desired.
\end{proof}

Finally, roots are unique in $\pi_1(\Sigma)$.

\begin{lemma}\label{lem:roots-unique}
  For $\gamma, \delta$ non-trivial elements of $\pi_1(\Sigma)$ and $k,
  l \in \ZZ_{> 0}$, if $\gamma^k = \delta^l$ then there is an element
  $\eta$ so that $\gamma = \eta^{\lcm(k,l)/k}$ and $\delta =
  \eta^{\lcm(k,l)/l}$.
\end{lemma}

One can prove this, for instance, by taking a hyperbolic or flat metric
on~$\Sigma$ and taking the geodesic representative for the conjugacy
class of $\gamma^k = \delta^l$, which will multiply-cover a
well-defined primitive loop~$L$.  Then $\gamma$ and $\delta$ must
both be conjugate to powers of the monodromy of~$L$.


\section{Skein algebra}
\label{sec:skein}

\begin{figure}
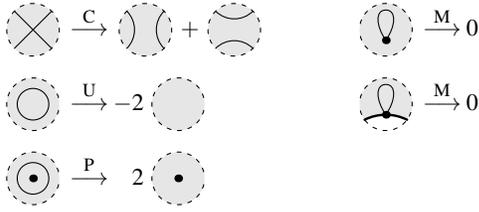

  \vskip-\abovedisplayskip
  \begin{align*}
  \mfigb{curves-0} &\move{C} \mfigb{curves-2} + \mfigb{curves-1}&
    \mfigb{curves-48} &\move{M} 0\\
  \mfigb{curves-10} &\move{U} -2\,\mfigb{curves-12}&
    \mfigb{curves-42} &\move{M} 0\\
  \mfigb{curves-11} &\move{P} \hphantom{-}2\,\mfigb{curves-13}
  \end{align*}
  \caption{The relations in $\Sk(\Sigma)$, when $\Sigma$ is a surface
    with punctures.  The relation in each case is that the left hand
    side equals the right hand side, but we draw an arrow to indicate
    that it is a reduction as we go from left to right.  For reduction
    (M) there can be any number of other
    arcs ending at the pictured marked point.
  \label{fig:skein-1}}
\end{figure}

The skein algebra $\Sk(\Sigma)$ is the quotient of $\ZZ\Diag(\Sigma)$ by the
relations from
Fig.~\ref{fig:skein-1}.  More precisely, we make the following
definitions.
\begin{definition}
  For $D \in \Diag(\Sigma)$, a \emph{skein reduction}
  of~$D$ is the element of $\ZZ\Diag(\Sigma)$ obtained by one of the
  replacements illustrated in Figure~\ref{fig:skein-1}:
  \begin{description}
  \item[Crossing resolution (C):] Replace a crossing by
    the sum of its two resolutions.
  \item[Unknot removal (U):] Replace an embedded $0$-gon by
    $-2$ times the rest of~$D$.
  \item[Punctured disk removal (P):] Replace a simple loop surrounding
    a puncture by $2$ times the rest of~$D$.
  \item[Monogon removal (M):] Replace an embedded monogon at a marked
    point by~$0$.
  \end{description}
  In each case, the \emph{reduction disk} is the
  disk indicated in Fig.~\ref{fig:skein-1}.  For reductions (C)\@, (U)\@, and
  (P) the intersection of~$D$ with the reduction disk is as shown;
  for (M) there may be other arcs ending at the marked point.

  Similarly, if $D$ skein reduces to $z$ and $z_3\in\ZZ\Diag(\Sigma)$
  does not have a term involving $D$, we also say that
  $\langle D \rangle
  + z_3$ has an elementary skein reduction to $z + z_3$.  Say
  that $z_1$ \emph{skein reduces} to $z_2$ if they differ by a
  sequence of zero or more elementary skein reductions, always going
  from left to right.
\end{definition}

\begin{definition}
  The \emph{skein algebra} $\Sk(\Sigma)$ is the quotient of
  $\ZZ\Diag(\Sigma)$ by skein reductions.  Let $\langle
  D \rangle$ be the image of a diagram~$D$ in $\Skein(\Sigma)$.
\end{definition}

\begin{proposition}\label{prop:isotopy-invariance}
  If $D_1$ and $D_2$ are regular isotopic, then $\langle D_1 \rangle =
  \langle D_2 \rangle$.
  If $D_1$ differs from $D_2$ by an RI move,
  then $\langle D_1 \rangle = -\langle D_2 \rangle$.
\end{proposition}

\begin{proof}
  Expand both sides using skein reductions.
\end{proof}

Prop.~\ref{prop:isotopy-invariance} lets us talk about the skein class
$\langle C \rangle$ of a curve~$C$.

\begin{proposition}\label{prop:product}
  On a punctured surface~$\Sigma$, union of curves induces a
  commutative,
  associative product on $\Sk(\Sigma)$.
\end{proposition}

\begin{proof}
  As noted in Lemma~\ref{lem:curves-prod}, there is a product on
  $\Curve(\Sigma)$.  By Prop.~\ref{prop:isotopy-invariance},
  this gives a map $\ZZ\Curve(\Sigma) \times \ZZ\Curve(\Sigma) \to
  \Sk(\Sigma)$.  We must show that for $C_1, C_1' \in \ZZ\Curve(\Sigma)$,
 if $\langle C_1\rangle = 
  \langle C_1'\rangle$ then $\langle C_1 C_2 \rangle = \langle C_1' C_2
  \rangle$.  It suffices to consider the case that $C_1$ is a
  single diagram (mod regular isotopy) and $C_1$
  and~$C_1'$ differ by an elementary reduction.  For reductions (C)\@,
  (U)\@, and~(M)\@, we may assume by ambient isotopy of~$C_2$ that the
  reduction disk in $C_1$ does not intersect~$C_2$, so the reduction descends.
  For reduction~(P)\@, we cannot avoid arcs of $C_2$ that end at the
  enclosed puncture, but a direct computation shows that
  the two sides have the same value.
\end{proof}

\begin{definition}\label{def:complexity}
  The \emph{complexity} of a diagram~$D$ is the pair
  \[
  (c(D), r(D))
  \]
  where
  $c(D)$ is the number of crossings of~$D$ and
  $r(D)$ is the number of reducible components of~$D$.
  A \emph{reducible
    component} of~$D$ is a place where moves (U)\@, (P)\@, or (M) can
  be applied.

  These pairs are ordered with the lexicographic order.  The
  complexity of a linear combination of diagrams is the list of
  complexities of the non-zero terms, sorted in decreasing order,
  and considered with the lexicographic order on lists.

  A diagram is \emph{irreducible} if its complexity is~$(0,0)$.
  Irreducible diagrams
  are nearly the same as simple diagrams, except that loops around
  punctures are forbidden.
\end{definition}

\begin{lemma}\label{lem:simplify}
  Complexity strictly reduces under skein reductions.
\end{lemma}

\begin{proof}
  This is immediate for reductions of a single diagram, and follows
  easily for linear combinations.
\end{proof}

\begin{lemma}\label{lem:chains-term}
  Any sequence of skein reductions in $\ZZ\Diag(\Sigma)$
  eventually
  terminates.
\end{lemma}

\begin{proof}
  The ordering on complexities
  is a well-ordering.
\end{proof}

\begin{proposition}\label{prop:simple-basis}
  Irreducible diagrams form a basis for $\Sk(\Sigma)$.  Any element of
  $\Sk(\Sigma)$ can be expressed in this basis by applying skein
  reductions whenever possible, in any order.
\end{proposition}

\begin{proof}
  Lemma~\ref{lem:chains-term} implies that irreducible
  curves span $\Sk(\Sigma)$.  To see linear independence, observe that
  reductions satisfy a \emph{diamond property}: if $z_1, z_2, z_3
  \in \ZZ\Diag(\Sigma)$ are such that $z_1$ has an elementary skein
  reduction to both $z_2$ and $z_3$, then there is another element
  $z_4 \in \ZZ\Diag(\Sigma)$ so that both $z_2$ and $z_3$ skein reduce
  to~$z_4$.  (This is very easy for this skein algebra, since the
  reducing disks for different relations
  can almost never overlap, except for monogons at the same puncture.)
  Then the diamond lemma \cite{Newman42:Equiv} implies that any two
  sequences of skein reductions terminate at the same place.
\end{proof}

\begin{remark}
  The reduction~(P) may look unfamiliar to readers used to the Jones
  skein.  This value for the loop is forced by consistency with
  multiplication (Prop.~\ref{prop:product}).  The quantum analogue is
  not obvious, and requires introducing \emph{opened surfaces} (cf.\
  \cite{FT12:ClusterII}).  Details will appear in a future paper.
\end{remark}


\section{Three bases}
\label{sec:three-bases}

Prop.~\ref{prop:simple-basis} gives a basis for the skein
algebra, but this basis is not always positive in the sense of
Def.~\ref{def:positive}.  There are, in fact, three related
bases.  We first give
the elementary building blocks.

\begin{definition}\label{def:basis-ops}
  For a simple loop~$L\in\Diag(\Sigma)$ with holonomy~$\gamma\in\pi_1(\Sigma)$ and an integer~$k>0$, we
  define three ways to create an element of $\QQ\Diag(\Sigma)$,
  as in
  Fig.~\ref{fig:three-bases}.
  \begin{itemize}
  \item $\Bang^k(L)$, the $k$'th \emph{bangle} of~$L$, is $L^k$, i.e., $k$
    parallel copies of~$L$.
  \item Let $I \subset L$ be a short interval.  Then
     $\Band^k(L)$, the $k$'th \emph{band} of~$L$, is $k$ copies
    of $L \setminus I$ with the
    ends connected by averaging
    all $k!$ ways of pairing the endpoints on the two sides.
  \item $\Brac^k(L)$, the $k$'th \emph{bracelet} of~$L$, is the
    loop whose holonomy is~$\gamma^k$, embedded tautly.  
    By Prop.~\ref{prop:powers-intersect}, $\Brac^k(L)$
    has $k-1$ self-intersections.
  \end{itemize}
  For convenience, for a simple loop~$L$ and simple arc~$A$, also define
  \begin{align*}
    \Bang^0(L) &= 1 & \Band^0(L) &= 1 & \Brac^0(L) &= 2\\
    \Bang^k(A) &= \langle A\rangle^k &
      \Band^k(A) &= \langle A\rangle^k &
        \Brac^k(A) &= \langle A \rangle^k.
  \end{align*}
\end{definition}
A priori, $\Band^k(L)$ is only in $\Sk_\QQ(\Sigma) \coloneqq
\Sk(\Sigma) \otimes \QQ$. In fact, it is in
$\Sk(\Sigma)$,
i.e., the coefficients are integral after reducing modulo the
skein relations (Prop.~\ref{prop:chebyshev-band}).

\begin{example}\label{examp:k=2}
  If
  $\langle L \rangle = z$, then
  \begin{align*}
    \Bang^2(L) &= z^2\\
    \Band^2(L) &= z^2 - 1\\
    \Brac^2(L) &= z^2 - 2.
  \end{align*}
  The first equation is trivial.  The third equation follows by applying the
  skein relation at the unique crossing of $\Brac^2(L)$.  The
  second equation is the average of the other two.
\end{example}

\begin{definition}\label{def:chebyshev}
  The \emph{Chebyshev polynomials of the first kind} are polynomials
  $T_n(z)$ satisfying the recurrence
  \begin{align}
    T_0(z) &= 2\label{eq:cheb1-0}\\
    T_1(z) &= z\label{eq:cheb1-1}\\
    T_{n+1}(z) &= z T_n(z) - T_{n-1}(z).\label{eq:cheb1-recur}
  \end{align}
  They satisfy
  \begin{align}
    T_k(z) T_l(z) &= T_{k+l}(z) + T_{\abs{k-l}}(z)\label{eq:cheb1-prod}\\
    T_k(e^x + e^{-x}) &= e^{kx} + e^{-kx}.
  \end{align}
  The \emph{Chebyshev polynomials of the second kind} are polynomials
  $U_n(z)$ satisfying the recurrence
  \begin{align}
      U_0(z) &= 1\label{eq:cheb2-0}\\
      U_1(z) &= z\label{eq:cheb2-1}\\
      U_{n+1}(z) &= zU_n(z) - U_{n-1}(z).\label{eq:cheb2-recur}
  \end{align}
 They satisfy
  \begin{align}
    U_k(z) U_l(z) &= U_{k+l}(z) + U_{k+l-2}(z) + \dots + U_{\abs{k-l}}(z)\label{eq:cheb2-prod}\\
    U_k(e^x + e^{-x}) &= e^{kx} + e^{(k-2)x} + \dots + e^{-kx}
  \end{align}
\end{definition}

\begin{proposition}\label{prop:chebyshev-bracelet}
  For any simple loop~$L$ and integer $n > 0$,
  \begin{equation}
    \Brac^n(L) = T_n(\langle L\rangle).
  \end{equation}
\end{proposition}

\begin{proof}
  This is trivially true for $n=1$ and the $n=2$ case was done in
  Example~\ref{examp:k=2}.
  To compute $\Brac^n(L)$ with $n > 2$, resolve one of the two outer
  crossings of $\Brac^n(L)$ (with holonomies $\gamma$ and
  $\gamma^{n-1}$, as in Prop.~\ref{prop:powers-intersect}).  One of the two resolutions is $L \cdot
  \Brac^{n-1}(L)$.  The other resolution differs by an RI move from
  $\Brac^{n-2}(L)$.  This gives
  Eq.~\ref{eq:cheb1-recur} in $\Sk(\Sigma)$.
\end{proof}

\begin{remark}
  By resolving $\Brac^n(L)$ at other crossings, we get a short proof
  of Eq.~\ref{eq:cheb1-prod} as applied to bracelets.
\end{remark}

To find the analogue of Prop.~\ref{prop:chebyshev-bracelet} for
the bands basis,
we introduce the graphical notation that
a box with an $n$ inside means averaging over all ways of joining
the $n$ strands on the
two sides of the box.
Diagrams should be interpreted as having a variable
number of strands (including~$0$), as indicated in the boxes.

\begin{lemma}\label{lem:box-identities}
  The following identities hold in $\Sk_\QQ(\Sigma)$:
  \begin{align}
    \label{eq:idempotent}
    \left\langle\mfigb{proj-10}\right\rangle &=
      \left\langle\mfigb{proj-0}\right\rangle\\
    \label{eq:reduce}
    \left\langle\mfigb{proj-0}\right\rangle
      &=  \frac{1}{n}\left\langle\mfigb{proj-2}\right\rangle
        + \frac{n-1}{n}\left\langle\mfigb{proj-3}\right\rangle\\
    \label{eq:reduce2}
      &= \left\langle\mfigb{proj-2}\right\rangle
        + \frac{n-1}{n}\left\langle\mfigb{proj-4}\right\rangle\\
    \label{eq:ptrace}
    \left\langle\mfigb{proj-5}\right\rangle
      &= -\frac{n+1}{n}\left\langle\mfigb{proj-6}\right\rangle.
  \end{align}
\end{lemma}

\begin{proof}
  Eq.~\ref{eq:idempotent} is true because averaging twice is
  the same as averaging once.  To see Eq.~\ref{eq:reduce}, note
  that if we average over $S_n$, the first strand on the top is
  connected to the first strand on the bottom with probability $1/n$,
  and connected somewhere else with probability $(n-1)/n$.  These two
  possibilities correspond to the two terms on the right.  Applying
  the skein relation at the crossing gives
  Eq~\ref{eq:reduce2}.  Eq.~\ref{eq:ptrace} follows from
  Eq.~\ref{eq:reduce2} by taking a partial trace: join the first
  strand on the top to the first strand on the bottom.
\end{proof}

\begin{remark}
  Lemma~\ref{lem:box-identities} is the $q=1$
  specialization of standard equations for the Jones-Wenzl
  idempotents \cite{Wenzl87:SeqProjections}.
\end{remark}

\begin{proposition}\label{prop:chebyshev-band}
  For any simple loop~$L$ and integer $n > 0$,
  \begin{equation}
    \Band^n(L) = U_n(\langle L\rangle)
  \end{equation}
\end{proposition}

\begin{proof}
  This is trivial for $n=1$ and already checked for $n=2$.
  For $n >2$, by
  Eqs.~\ref{eq:reduce2}, \ref{eq:idempotent},
  and~\ref{eq:ptrace} (in that order) to the diagram defining
  $\Band^n(L)$, we find that $\Band^n(L)$
  satisfies Eq.~\ref{eq:cheb2-recur}.
\end{proof}

\begin{definition}\label{def:bases}
  The \emph{bangles basis} $\Basis_0(\Sigma)$ for $\Sk(\Sigma)$
  consists of all irreducible diagrams.  It is parameterized by \emph{integer
    laminations}, which are unordered collections
  $\mu = \{(a_i, C_i)\}_{i=1}^k$, where
  \begin{itemize}
  \item each $a_i$ is a positive integer,
  \item each $C_i$ is a connected irreducible diagram,
  \item no two $C_i$ intersect, and
  \item all the $C_i$ are distinct up to ambient isotopy.
  \end{itemize}
  The corresponding basis element in $\Basis_0(\Sigma)$ is
  \[
  \Bang(\mu) \coloneqq \prod_{i=1}^k \Bang^{a_i}(C_i) = \prod_{i=1}^k
  \langle C_i\rangle^{a_i}.
  \]
  The \emph{bands} and \emph{bracelets bases} $\Basis_1(\Sigma)$ and
  $\Basis_2(\Sigma)$ are also parameterized by
  integer laminations, but with corresponding basis elements
  \[
  \Band(\mu) \coloneqq \prod_{i=1}^k \Band^{a_i}(C_i)\qquad
  \Brac(\mu) \coloneqq \prod_{i=1}^k \Brac^{a_i}(C_i),
  \]
  respectively.  Informally, start from an irreducible diagram and
  replace parallel loops with the corresponding band
  or bracelet.
\end{definition}

\begin{proposition}
  The sets $\Basis_0(\Sigma)$, $\Basis_1(\Sigma)$, and
  $\Basis_2(\Sigma)$ each form a basis for $\Sk(\Sigma)$.
\end{proposition}

\begin{proof}
  For $\Basis_0(\Sigma)$ this is
  Prop.~\ref{prop:simple-basis}.  The other two bases are
  triangular with respect to $\Basis_0(\Sigma)$.
\end{proof}


Here is an intrinsic characterization of the bracelets basis.

\begin{definition}
  A \emph{bracelet loop} is a loop of the form $\Brac^k(L)$ for some
  $k \ge 1$ and some simple loop~$L$.
  A \emph{multi-bracelet} is a diagram~$D$ in which
  \begin{itemize}
  \item each component is a simple arc or a bracelet loop, and
  \item no two components intersect.
  \end{itemize}
  A multi-bracelet has \emph{parallel bracelets} if there are two
  components that are bracelets of the same simple loop.
\end{definition}

\begin{lemma}
  A diagram is in $\Basis_2(\Sigma)$ iff it is a multi-bracelet
  with no parallel bracelets and no bracelets of punctured disks.
\end{lemma}

\subsection{Examples on positivity}
Theorem~\ref{thm:positive-1} says
that the bracelets basis is positive.  We will now give some other
examples and conjectures on
positivity and
non-positivity.

\begin{example}
  Let $\Sigma$ be the annulus with two marked points, one on each
  boundary component.  Let $A_k$ be the arc connecting the two marked
  points at slope~$k$ (i.e., wrapping $\abs{k}$ times around the core
  loop, clockwise or counterclockwise according to the sign of~$k$).
  Let $L$ be the core loop, and let $B$ be a push-off of the union of the two
  boundary components (considered as a diagram with two arcs).
  Then elementary induction using the skein rules shows that
  \begin{align}
    \langle\Brac^n(L)\rangle \cdot \langle A_k\rangle &=
      \langle A_{k+n}\rangle + \langle A_{k-n}\rangle\\
    \langle A_a\rangle \cdot \langle A_{a+k}\rangle &=
        \langle B \rangle \sum\limits_{i=1}^{\lfloor k/2\rfloor} i\Brac^{k-2i}(L)\\
        &\qquad\qquad
      + \begin{cases}
        \bigl\langle A_{a+\frac{k}{2}}\bigr\rangle^2
        & \text{$k$ even}\\
        \bigl\langle A_{a+\frac{k-1}{2}}\bigr\rangle\bigl\langle
        A_{a+\frac{k+1}{2}}\bigr\rangle
        & \text{$k$ odd}.
      \end{cases}\nonumber
  \end{align}
  The structure constants are positive, as expected.  They
  are also positive in the bands basis, but in the bangles basis,
  \begin{align*}
    \langle A_0 \rangle \langle A_5 \rangle
    &= \langle A_2 \rangle \langle A_3\rangle +
    \langle B \rangle \bigl(\langle L\rangle^3 -
    \langle L\rangle\bigr).
  \end{align*}
\end{example}

\begin{example}\label{examp:torus}
  Take $\Sigma = T^2$, the unpunctured torus.
  Let $\alpha, \beta$ be the two generators for $\pi_1(T^2)$, and for
  $\gamma \in \pi_1(T^2)$, $\gamma \ne 1$, let $\gamma^\circ\in\Diag(\Sigma)$ be a taut
  representative for the corresponding conjugacy class.
  Set $C_{a,b} = (\alpha^a \beta^b)^\circ$ for $a,b \in \ZZ$, not
  both~$0$.  The bracelets basis is
  \begin{align*}
    \Basis_2(T^2) &= \{\,\langle C_{a,b}\rangle \mid (a,b) \ne (0,0)\,\} \cup
      \{\,1\,\}.
  \end{align*}
  The only duplicates on the list above arise from the equality
  $\langle C_{a,b}\rangle = \langle C_{-a,-b}\rangle$.
  For convenience, also define $\langle C_{0,0}\rangle = 2$.

  \begin{proposition}\label{prop:torus-mult}
    The multiplication rules for $\Sk(T^2)$ are
    \begin{equation}
    \langle C_{a,b}\rangle \cdot \langle C_{c,d}\rangle =
     \langle C_{a+c,b+d}\rangle + \langle C_{a-c,b-d}\rangle.
    \label{eq:torus-mult}
    \end{equation}
  \end{proposition}

  \begin{proof}
    This follows from \cite{FG00:SkeinTorus}, but we give a short
    independent argument.
    If $(a,b) = \pm (c,d)$, Eq.~\ref{eq:q-torus-mult}
    follows from Prop.~\ref{prop:chebyshev-bracelet}.
    Otherwise, take a taut embedding of
    $(\alpha^a\beta^b)^\circ
    \cup (\alpha^c\beta^d)^\circ$ and resolve any one crossing between the
    two components.  Both resolutions
    are strongly positive (see
    Lemma~\ref{lem:distinct-simple-positive} below) and are $(\alpha^{a+c}\beta^{b+d})^\circ$ and
    $(\alpha^{a-c}\beta^{b-d})^\circ$.  
  \end{proof}

  \begin{corollary}
    The bands basis and the bangles basis are not positive on $T^2$.
  \end{corollary}

  \begin{proof}
    Prop.~\ref{prop:torus-mult} tells us that
    \begin{align*}
       \langle C_{0,1} \rangle  \langle C_{2,1} \rangle &=
        \Brac^2(C_{1,0}) +  \Brac^2(C_{1,1})\\
      &= \Band^2(C_{1,0}) + \Band^2(C_{1,1}) - 2\\
      &= \Bang^2(C_{1,0}) + \Bang^2(C_{1,1}) - 4. \qedhere
    \end{align*}
  \end{proof}
\end{example}

\begin{example}\label{examp:perf-torus}
  We can do some computations for $\Sigma=T^2 \setminus D^2$, a
  perforated torus.
  Let $B$ be the loop around the boundary.  For $(k,l) \in
  \ZZ^2$ relatively prime, let $C_{k,l}$ be the simple loop at slope
  $k/l$ (i.e., homotopic to $\alpha^k\beta^l$)
  on the unpunctured torus.  (There is a unique such loop.)  The
  first interesting product is
  \begin{align*}
    \langle C_{0,1}\rangle \langle C_{2,1}\rangle
    &= \langle C_{1,1}\rangle^2 + \langle C_{0,1} \rangle^2 + \langle B\rangle - 2\\
    &= \Band^2(C_{1,1}) + \Band^2(C_{0,1}) + \langle B\rangle\\
    &= \Brac^2(C_{1,1}) + \Brac^2(C_{0,1}) + \langle B\rangle+2.
  \end{align*}
  This is positive in the bracelets and bands bases, but not for the
  bangles basis, and reduces to the answer from
  Eq.~\ref{eq:torus-mult} if we set $\langle B \rangle = -2$.
\end{example}

\begin{problem}
  Find a general formula for $\langle C_{a,b}\rangle \langle
  C_{c,d}\rangle$ on $T^2 \setminus D^2$.
\end{problem}

The examples above imply that the bangles basis is almost never
positive.  On the other hand, the following conjecture is
plausible.

\begin{conjecture}\label{conj:bands-positive}
  The bands basis is positive when $\Sigma$ has either
  non-empty boundary or at least one puncture, i.e., if
  $\pi_1(\Sigma)$ is free.
\end{conjecture}

Conjecture~\ref{conj:bands-positive} includes all cases related to
cluster algebras.

\subsection{Non-commutative skein algebra}
We use an extension of the
usual non-commutative Jones skein algebra to allow
marked points on
$\partial\Surf$ as described by Muller~\cite{Muller12:SkeinCluster}.
(There is a further extension that allows
punctures, which will be the subject of a future paper.)
In that setting, $\Sk_q(\Sigma)$
also has a basis consisting of simple
diagrams.  The three bases generalize: for a simple loop~$L$, define
in $\Sk_q(\Sigma)$
\begin{align*}
  \Bang^k(L) &= \langle L\rangle^k\\
  \Band^k(L) &= U_k(\langle L\rangle)\\
  \Brac^k(L) &= T_k(\langle L\rangle).
\end{align*}
Extend this to a complete basis as before.
Note that $U_k$ and $T_k$ are the ordinary Chebyshev
polynomials with integer coefficients.

We say that a basis for an algebra over $\ZZ[q^\pm]$ is
\emph{positive} if the structure constants for multiplication lie in
$\ZZ_{\ge 0}[q^\pm]$.

\begin{conjecture}\label{conj:bracelets-q-positive}
  The bracelets basis is a
  positive basis for $\Sk_q(\Sigma)$.
\end{conjecture}

Evidence for Conjecture~\ref{conj:bracelets-q-positive} comes from
$\Sk_q(T^2)$, where Frohman and Gelca computed the multiplication
rules with respect to a
suitable basis \cite{FG00:SkeinTorus}:
  \begin{equation}
    \label{eq:q-torus-mult}
    C_{a,b} \cdot C_{c,d} =
     q^{ad-bc} C_{a+c,b+d} +
     q^{-ad+bc} C_{a-c,b-d}.
  \end{equation}

\begin{conjecture}\label{conj:bands-q-positive}
  For any connected surface $\Sigma$ with non-empty boundary, the
  bands basis is a positive basis
  for $\Sk_q(\Sigma)$.
\end{conjecture}

\begin{remark}
  Although the bracelets basis is the subject of this
  paper, the
  bands basis has appeared in several different contexts, and arises
  naturally from representation theory.  When we
  think of $\Sk(\Sigma)$ as functions on the twisted $\SL_2(\RR)$
  character variety of~$\Sigma$, the value of a single loop~$L$ with
  lifted holonomy $\gamma \in \pi_1(UT\Sigma)$ is
  \[
  \langle L \rangle = \pm \tr_2(\rho(\gamma))
  \]
  where $\rho$ is the $\SL_2$ representation and $\tr_2$ is the trace
  in the $2$-dimensional representation of
  $\SL_2$.  With this setup, we
  have
  \begin{align*}
    \Band^k(L) &= \pm\tr_{k+1}(\rho(\gamma))\\
    \Brac^k(L) &= \pm\tr_2(\rho(\gamma^k)).
  \end{align*}
  That is, for $\Band^k(L)$, we take the trace in the $k$'th symmetric
  power of the defining representation of $\SL_2$ (with dimension
  $k+1$), while for $\Brac^k(L)$ we take the trace of the $k$'th power
  of the loop.

  In the case of an annulus with two marked points (one on each
  boundary component), the skein algebra is contained in a quantized
  affine algebra, and Lusztig's dual canonical basis is the bands
  basis and not the bracelets basis~\cite{Lampe11:Kronecker}.
\end{remark}


\section{Positivity}
\label{sec:positive}

\begin{definition}
  A diagram~$D$ is (manifestly) \emph{null} if it has a singular
  monogon based at a marked point.
  A crossing $x$ of a diagram~$D$ is
  \emph{weakly positive} if both resolutions of~$D$ at~$x$
  are weakly positive or null, with at
  least one of the resolutions weakly positive; similarly for strongly
  positive.
\end{definition}

Here is the key lemma of the paper.

\begin{lemma}\label{lem:key}
  Let $D$ be a weakly positive diagram that is not isotopic to a
  multi-bracelet.  Then $D$ has weakly positive crossing.
\end{lemma}

\begin{proof}[Proof of Theorem~\ref{thm:positive-expansion}, assuming Lemma~\ref{lem:key}]
  First suppose we have a weakly positive diagram~$D$.  We proceed by
  induction on the complexity of~$D$
  (see Def.~\ref{def:complexity}).

  If $D$ is not isotopic to a multi-bracelet, apply Lemma~\ref{lem:key}:
  since $D$ has a weakly positive crossing, we have
  $\langle D\rangle = \langle D_1\rangle + \langle D_2 \rangle$
  where both $D_i$ are weakly positive or null and strictly simpler
  than~$D$.  Null diagrams can be ignored (as they are $0$ in the
  skein algebra).
  So, by induction, the $\langle D_i\rangle$ and therefore
  $\langle D \rangle$ have positive
  expansions in the bracelets basis.

  If $D$ is isotopic to a multi-bracelet~$D'$, it is nearly in the
  bracelets basis,
  except that it may have parallel bracelets or bracelets around punctures.
  Any bracelet around a puncture is equal to~$2$, so these may be removed.
  If $D'$ has parallel bracelets, say $D' = \Brac^k(L)\Brac^l(L) D''$ for
  a simple loop~$L$, then by
  Prop.~\ref{prop:chebyshev-bracelet} and
  Eq.~\ref{eq:cheb1-prod},
  \begin{equation}
    \langle D \rangle = \langle D' \rangle
      = \langle D'' \rangle \langle \Brac^{k+l}(L) \rangle
        + \langle D'' \rangle \langle \Brac^{\abs{k-l}}(L) \rangle.
  \end{equation}
  The terms on the right have fewer parallel bracelets than $D'$, and so we
  can repeat this reduction until we are left with a positive linear
  combination of
  elements of the bracelets basis.

  If $D$ is not positive, then by Theorem~\ref{thm:curve-shorten} it
  can be reduced to a taut curve $D'$ by a sequence of Reidemeister
  moves.  Each move $D_i \to D_{i+1}$ either leaves the value in the skein algebra the
  same or multiplies it by an integer, as follows:
  \begin{description}
  \item[Move RO:] $\langle D_i\rangle = -2\langle D_{i+1}\rangle$.
  \item[Move RI:] $\langle D_i\rangle = -\langle D_{i+1}\rangle$.
  \item[Move RIb:] $\langle D_i\rangle = 0$.
  \item[Moves RII\@, RIIb, and RIII\@:] $\langle D_i \rangle = \langle D_{i+1}\rangle$.
  \end{description}
  Thus $\langle D\rangle = k\langle D'\rangle$ for some integer~$k$, and the theorem follows.
\end{proof}

\begin{proof}[Proof of Theorem~\ref{thm:positive-1}, assuming Lemma~\ref{lem:key}]
  Given two diagrams $D_1$ and~$D_2$ in the bracelets basis, find a taut
  representative $D$ for $D_1 \cup D_2$ and apply
  Theorem~\ref{thm:positive-expansion}.
  Lemma~\ref{lem:taut-union} guarantees that $\langle D\rangle =
  \langle D_1\rangle\langle D_2\rangle$.
\end{proof}

So we only need to prove Lemma~\ref{lem:key}.  We prove a slightly
stronger version.

\begin{lemma}\label{lem:key2}
  Let $D$ be a taut diagram that is not a
  multi-bracelet.  Then $D$ is isotopic through RIII moves
  to a taut diagram~$D'$ that has a strongly positive crossing.
\end{lemma}

\begin{proof}[Proof of Lemma~\ref{lem:key}, assuming Lemma~\ref{lem:key2}]
  Let $D$ be a weakly positive diagram.  Then by assumption $D$ can
  be reduced by
  regular isotopy to a taut diagram~$D'$, which by
  Lemma~\ref{lem:key2} is regular isotopic to a diagram~$D''$
  with a strongly positive crossing~$x$.  By
  Lemma~\ref{lem:strong-weak-positive}, $x$ is also weakly positive.  By
  Lemma~\ref{lem:cross-correspond}, there is a crossing
  $m(x)\in\Cross(D)$ so that
  the resolutions of $m(x)$ are regular isotopic to the resolutions
  of~$x$.  Thus $m(x)$ is also weakly positive.
\end{proof}

The plan to prove Lemma~\ref{lem:key2} is to
look for a crossing at the end of a \emph{maximal bracelet chain}, a
portion of the diagram that wraps around a simple loop
(Def.~\ref{def:bracelet} below), like
crossing $x_3$ in Fig.~\ref{fig:partial-bracelet}.
This crossing will be strongly positive
(Lem.~\ref{lem:max-brac-positive}).  If there are no such
crossings, each
component of the diagram is a simple arc or a bracelet and any
crossing between two components is strongly positive
(Lem.~\ref{lem:distinct-simple-positive}).

First we build up some tools.

\begin{lemma}\label{lem:reduction-cond}
  Suppose $D$ is a taut diagram and $D \to D_1 + D_2$ is a crossing reduction
  with $D_1$ not strongly positive.  Then there is a $0$-gon
  or $1$-gon $H \subset D_i$ passing through the
  reduction disk twice.
\end{lemma}

\begin{proof}
  By definition, if $D_1$ is not strongly positive it has a singular disk or
  monogon~$H$.
  If $H$ does not pass through the reduction disk, then
  it is also a singular monogon or disk for~$D$, contradicting the
  assumption that $D$ is taut.
  If $H$ passes
  through the reduction disk once, then we can lift it to~$D$:
  a $0$-gon lifts to a $1$-gon in~$D$, and a $1$-gon lifts to a $2$-gon in~$D$.
  But a taut diagram cannot have any $1$-gons, and
  by Lemma~\ref{lem:no-gons} can only have
  $2$-gons between two parallel strands.
  So $H$ passes through the reduction disk at least twice.  Since by
  definition the intervals making up a chain are disjoint, $H$
  must pass through the reduction disk exactly twice.
\end{proof}

\begin{lemma}\label{lem:discon-positive}
  Suppose $D$ is a taut diagram, $x \in \Cross(D)$, and $D_1$ is the
  disconnected resolution of $D$ at~$x$ (i.e., the resolution that
  increases the number of components).  Then $D_1$ is strongly positive.
\end{lemma}

\begin{proof}
  By Lemma~\ref{lem:reduction-cond}, if $D_1$ is not strongly positive
  there is a $0$-gon or $1$-gon passing through the reduction disk
  twice.  But a $0$-gon or $1$-gon cannot touch both components of the
  disconnected resolution, so this is impossible.
\end{proof}

\begin{definition}\label{def:bracelet}
  In a taut diagram~$D$, recall that a \emph{bracelet loop} is a
  component of~$D$ that is homotopic to $\Brac^k(L)$ for a simple loop~$L$.  A
  \emph{bracelet chain} is a $1$-chain~$H$, in the sense of
  Def.~\ref{def:chains-gons}, so that $H^\circ$ is homotopic to
  $\Brac^k(L)$ for some simple loop~$L$ and $k > 0$ and with $H$ not a
  complete arc component of~$D$.  A \emph{bracelet} is a bracelet
  loop or bracelet chain.  A bracelet~$B$ has
  an \emph{underlying loop}~$L(B)$, around which it travels~$n(B)$ times.
  Define
  $\gamma(B)\in\pi_1(\Sigma)$ to be the holonomy of $L(B)$, with
  basepoint and arc to the basepoint specified as necessary. The
  holonomy of $B$ itself is $\gamma^n(B)\coloneqq \gamma(B)^{n(B)}$.  Finally, a
  \emph{maximal bracelet} is a bracelet~$B$ for which there is no
  other bracelet~$B'$ with $L(B')=L(B)$ and $n(B') > n(B)$.
\end{definition}

\begin{figure}
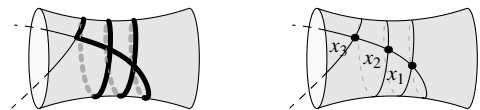

  \vskip-\abovedisplayskip
  \begin{align*}
    &\mfigb{annulus-12}
    &&\mfigb{annulus-10}
 \end{align*}
  \caption{A maximal bracelet chain and
    the labeling of the self-intersections.\label{fig:partial-bracelet}}
\end{figure}

\begin{lemma}\label{lem:max-bracelet}
  Let $D$ be a connected taut diagram with at least one
  self-intersection.  Then $D$ has a maximal bracelet.
\end{lemma}

\begin{proof}
  Among all $1$-chains~$H$ of $D$ with corners at crossings, pick one that
  is minimal with
  respect to inclusion.  (There is one since $D$ has
  a self-intersection.)  Then $H^\circ$ is necessarily a simple
  loop, so $H$ is a bracelet.  Now take
  $B$ to be a maximal bracelet with $L(B) = H^\circ$.
\end{proof}

We can
arrange for maximal bracelets to lie in a good position.

\begin{lemma}\label{lem:taut-bracelet}
  Let $D$ be a taut diagram and let $B \subset D$ be a bracelet chain.
  Then $(B,D)$ is isotopic through RIII moves to $(B',D')$
  where $D'$ and $(B')^\circ$ are both taut.
\end{lemma}

\begin{proof}
  Suppose that $B^\circ$ is not taut.  We will reduce the number of
  self-intersections of~$B^\circ$ by RIII moves on~$D$.

  Let $x$ be the intersection at the
  end of~$B$ and let $\bar x$ be a point on $B^\circ$ near~$x$.
  Note that $B^\circ$ is part of the disconnected resolution of $D$
  at~$x$, so by Lemma~\ref{lem:discon-positive} it is strongly
  positive.  Thus, by Theorem~\ref{thm:gon-exist-2}, $B^\circ$ has
  a singular $2$-gon~$G$.  Since $D$ is taut, $G$ must pass
  over~$\bar x$; let $G = (S_1, S_2)$, with the segment $S_1$
  containing $\bar x$.  Pick $G$ so that $S_1$ is minimal with respect
  to inclusion.

  Since $G$ is null-homotopic, it lifts to a $2$-gon $\t G = (\t S_1,
  \t S_2)$ in the universal cover $\t \Sigma$ of~$\Sigma$.  Since
  $S_1$ is minimal, $\t S_1$ and $\t S_2$ do not intersect, so $\t G$
  is an embedded 2-gon.  Let $\t A$ be the intersection of the preimage
  of $D$ with the disk bounded by $\t G$.  Because $D$ is taut, no
  component of~$\t A$ can meet $\t S_2$ twice, and because $S_1$ is
  minimal, no component of~$\t A$ can meet $\t S_1$ twice.  So each
  component of~$\t A$ is an arc running from $\t S_1$ to $\t S_2$.  Now
  let $D'$ be~$D$, but with $S_2$ pushed over~$G$ to run
  parallel to~$S_1$ outside of~$G$, and let $B'$ be the bracelet in $D'$
  with corner at~$x$.  Then $D'$ has the same number of
  self-intersections as~$D$, but $(B')^\circ$ has fewer
  self-intersections, as desired.

  The homotopy from $S_2$ to a parallel copy of $S_1$ can be done
  using only RIII moves, since it doesn't change the intersection number
  \cite[Lemma 1.6]{HS94:ShorteningCurves}.
\end{proof}

Following Lemma~\ref{lem:taut-bracelet}, we say that a bracelet chain $B
\subset D$ is \emph{taut} if $D$ and $B^\circ$ are both taut.  In this
case, we follow Prop.~\ref{prop:powers-intersect} and label the
self-intersection of~$B$ that cuts off holonomy $\gamma(B)^k$ by
$x_k(B)$ for $k=1,\dots,n$, as in Fig.~\ref{fig:partial-bracelet}.
(To define the holonomy, pick a basepoint on~$B$ and
connect the intersection to the basepoint by traveling along~$B$.)

\begin{lemma}\label{lem:max-brac-positive}
  Let $B \subset D$ be a taut maximal bracelet that is a chain.  Then
  for $k > n(B)/2$,
  the crossing $x_k(B)$ is
  strongly positive.
\end{lemma}

\begin{proof}
  For convenience, we assume that $D$ has only one connected
  component.
  Suppose $D'$ is a resolution of~$D$ at $x_k(B)$ that is
  not strongly positive.   By Lemma~\ref{lem:discon-positive}, $D'$ is the
  connected resolution.
  Let $B'$ be the sub-bracelet of~$B$ cut off by $x_k(B)$.
  Let $H$ be the $0$-gon from Lemma~\ref{lem:reduction-cond}.  If $H$
  is a 1-gon, its endpoints are either
  both in $B'$ or both in $D \setminus B'$.  We write $\gamma$ for
  $\gamma(B)$ and proceed by cases on~$H$.
  \begin{itemize}
  \item If $H$ is a $0$-gon, then $D$ is a loop with holonomy
    $\gamma^{2k}$.
  \item If $H$ is a $1$-gon with endpoints in $B'$, the endpoints are at $x_l(B)$
    for some $l < k$.  Let the holonomy of $D \setminus B'$ (as a $1$-chain)
    be~$\rho$.  Then the holonomy of $H$ is
    $\rho^{-1}\gamma^{k-l} = 1$, which implies that $D$ is a loop
    with holonomy
    $\rho\gamma^k = \gamma^{2k-l}$.
  \item If $H$ is a $1$-gon with endpoints in $D \setminus B'$, let
    $y$ be the corresponding corner.
    There is a corresponding $1$-chain $C$ in $D$ with corner at~$y$.
    Let $\rho$ be the holonomy of~$C\setminus B'$.
    Then the holonomy of $H^\circ$ is $\rho \cdot
    \gamma^{-k} = 1$ so $\rho = \gamma^k$, and the holonomy of $C$ is
    $\gamma^{2k}$.
  \end{itemize}
  In almost all cases, we found a bracelet chain which contradicts the
  maximality of~$B$.
  The only remaining case is when
  $k=(n+1)/2$ and $B$ is contained in an arc with holonomy conjugate
  to $\gamma^{n+1}$.  In this case, the connected
  resolution is null. The disconnected resolution is
  strongly positive by Lemma~\ref{lem:discon-positive}, so again the
  crossing is strongly positive.
\end{proof}

\begin{lemma}\label{lem:distinct-simple-positive}
  Let $D_1\cup D_2$ be a taut diagram where $D_1$ and $D_2$ are
  simple arcs or bracelet loops.  Then any
  crossing between $D_1$ and $D_2$ is strongly positive.
\end{lemma}

\begin{proof}
  Let $x$ be a crossing between $D_1$ and $D_2$, let $D'$ be a
  resolution of
  $D_1 \cup D_2$ at~$x$, and suppose $H$ is a $0$-gon or $1$-gon
  of~$D'$.  By Lemma~\ref{lem:reduction-cond}, $H$ must pass through
  the reduction disk at~$x$ twice, which means that one of the curves
  (say, $D_2$) must be a bracelet loop and $H$ must run completely
  over~$D_2$, with both endpoints in~$D_1$.

  If $D_1$ is an arc, then (since $D_1$ is simple) the endpoints
  of~$H$ must be the endpoints of~$D_1$, which are necessarily at
  the same marked
  point.  Let the holonomy of $(D_1)^\circ$ be $\rho\in\pi_1(\Sigma,x)$.
  If $\rho$ is a power of $\gamma(D_2)$, $D_1$ and $D_2$ do not
  intersect.  But otherwise the holonomy of~$H^\circ$
  is $\rho\cdot \gamma^{\pm n}(D_2)$, which is not~$1$ by
  Lem.~\ref{lem:roots-unique}.

  If $D_1$ and $D_2$ are both bracelet loops, again $\gamma(D_1) \ne
  \gamma(D_2)$ in $\pi_1(\Sigma, x)$, since otherwise $D_1$ and $D_2$
  would not intersect.
  The holonomy of $H^\circ$ is
  $\gamma(D_1)^k\cdot\gamma(D_2)^{n(D_2)}$ for some $1 \le k \le n(D_1)$.  Again
  this cannot be $1$ by Lemma~\ref{lem:roots-unique}.
\end{proof}

\begin{proof}[Proof of Lemma~\ref{lem:key2}]
  If any component~$C$ of~$D$ is not a simple arc or bracelet loop, then by
  Lemma~\ref{lem:max-bracelet} it has a maximal bracelet chain~$B$,
  which we can assume (after RIII moves) to be taut by
  Lemma~\ref{lem:taut-bracelet}.  Then $x_n(B)$ is
  strongly positive by Lemma~\ref{lem:max-brac-positive}.

  If all components of~$D$ are simple arcs or bracelets loops, by
  Lemma~\ref{lem:distinct-simple-positive}, any crossing between different
  components of~$D$ is strongly positive.

  Thus if $D$ is not a multi-bracelet, we have exhibited a strongly positive
  crossing in a diagram isotopic to~$D$ through RIII moves.
\end{proof}

\begin{remark}
  A closer analysis of the proof
  shows that every taut diagram~$D$ has a strongly positive crossing
  (i.e., we do not need to do an isotopy first): a crossing~$x$
  near the end of a
  maximal bracelet chain~$B$ is strongly positive, whether or not $B$
  is taut.  Essentially, the RIII moves to make~$B$ taut
  (following Lemma~\ref{lem:taut-bracelet}) are also RIII moves on the
  connected resolution of~$D$ at~$x$.
\end{remark}


\begin{acknowledgments}
I would like to thank
Thang Le,
Bruno Martelli,
Greg Muller,
Gregg Musiker, and
Lauren Williams for many helpful conversations.
This paper grew out of joint work with
Sergei Fomin and Michael Shapiro, who contributed greatly to the
development of this work.
The research for this project was largely carried out while visiting UC
Berkeley and attending the MSRI program on Cluster Algebras in Fall
2012, and I thank these institutions for their wonderful hospitality.
This work
was partially supported by NSF grant DMS-1008049.
\end{acknowledgments}

  \bibliographystyle{mypnas}
  \bibliography{curves,dylan,cluster,quantum,topo,misc,drafts}
\end{article}

\end{document}